\documentclass[12pt,reqno]{amsart}
\usepackage{mathrsfs}
\usepackage{amsfonts}
\usepackage{amsmath}
\numberwithin{equation}{section}
\newtheorem{theorem}{Theorem}[section]

\newtheorem{lemma}{Lemma}[section]
\newtheorem{proposition}{Proposition}[section]
\newtheorem{definition}{Definition}[section]
\newtheorem{remark}{Remark}[section]

\begin{document}
\title[Initial boundary value problem for  nonlinear  Dirac  equation]
      {Initial boundary value problem for nonlinear  Dirac  equation of Gross-Neveu type in $1+1$ dimensions}

\author{Yongqian Zhang}
\address{Yongqian Zhang:
School of Mathematical Sciences,\
\small Fudan University, Shanghai 200433, P.R.China}
\email{\tt yongqianz@fudan.edu.cn}

\author{Qin Zhao}
\address{Qin Zhao : School of Mathematical Sciences, \
\small Shanghai Jiao Tong University, Shanghai 200240, P.R.China}
\email{\tt zhao@sjtu.edu.cn}


\begin{abstract}
This paper studies an initial boundary value problem for a class of nonlinear  Dirac equations with cubic terms and moving boundary. For  the initial data with bounded $L^2$ norm and the suitable boundary  conditions, the global existence and the uniqueness of the strong solution are proved.
\end{abstract}

\maketitle

\section{Introduction}

Consider the nonlinear  Dirac equations
\begin{equation}\label{eq-dirac}
\left\{ \begin{array}{l} i(u_t+u_x)=-mv +N_1(u,v), \\i( v_t-v_x)=-mu+N_2(u,v),
\end{array}
\right.
\end{equation}
in a domain $\Omega=\{ (x,t)\, \big|\, t\ge 0, \, x\ge z(t)\}$ for $m\ge 0$
with initial data
\begin{equation}\label{eq-dirac-initialv}
(u(x,t=0), v(x,t=0))=(u_0(x), v_0(x)), \quad x\ge 0,
\end{equation}
and boundary condition
\begin{equation}\label{eq-dirac-bc}
u(z(t),t)=\lambda(t)v(z(t),t), \quad t\ge 0.
\end{equation}
The nonlinear terms  take the following form
\begin{eqnarray} N_1=\partial_{\overline{u}}W(u,v)=\alpha u|v|^2+ 2\beta v(\overline{u}v+u\overline{v}), \label{eq-nonlinearstruc1} \\ N_2=\partial_{\overline{v}}W(u,v)=\alpha v|u|^2+ 2\beta u(\overline{u}v+u\overline{v}) , \label{eq-nonlinearstruc2}
\end{eqnarray}
with
\[ W(u,v)=\alpha |u|^2|v|^2+\beta (\overline{u}v+u\overline{v})^2,\]
where $\alpha,\beta\in R^1$ and $\overline{u}, \overline{v}$ are complex conjugate of $u$ and $v$.

The boundary $\{x=z(t)\}$, denoted by $\Gamma_B$,  is assumed to satisfy the following,
\begin{description}
\item[(H1)] $-1<z_t(t)<1$, for $t\ge 0$ and $z(0)=0$.
\item[(H2)] $|\lambda(t)|^2(1-z_t(t))\le (1+z_t(t))$, for $t\ge 0$.
\end{description}
Here and in sequel, we denote $z_t=\frac{dz}{dt}$, $\lambda_t=\frac{d\lambda}{dt}$, $u_t=\frac{\partial u}{\partial t}$, $u_x=\frac{\partial u}{\partial x}$ etc. for simplification.

The nonlinear Dirac equation (\ref{eq-dirac}) is called  Thirring equation for $\alpha=1$ and $\beta=0$, while it is called Gross-Neveu equation for $\alpha=0$ and $\beta=1/4$; see for instance \cite{thirring} and \cite{gross-neveu}, \cite{pelinovsky}. There are a number of works devoted to the local and global well-poedness of the Cauchy problem for the nonlinear Dirac equation with various type of nonlinearities in different spatial dimensions (see for instance \cite{bachelot-2,bournaveas-zouraris,cacciafesta1,candy,deldado,dias-figueira2,escobedo,Esteban-Lewin-Sere, gross-neveu,huh2, pelinovsky,selberg,thirring,zhang,Zhang-zhao}, and the references therein). There are also some papers on the initial boundary value problem (see for example \cite{bouvier-gerad} and \cite{naumkin}). In \cite{bouvier-gerad}, motivated to study the Hawking effect describing the collapse of a spherically symmetric star to a Schwarzchild black hole, Bouvier and G\'{e}rard used  technique from $C^*$ algebra to study the asymptotic behaviour of the  global solution to  (\ref{eq-dirac}),(\ref{eq-dirac-initialv}) and (\ref{eq-dirac-bc}) with a class of special initial data in $R^{1+1}$, where the non-characteristic boundary is assumed to approach characteristic as $t\to\infty$, with $|\lambda(t)|^2(1-z_t(t))= (1+z_t(t))$ for $t\ge 0$ and the solution is assumed to be bounded.  In \cite{naumkin,naumkin16jde}, Naumkin proved the existence of global solution in $H^1$ to initial boundary value problem  for Thirring model in quarter plane $\{t>0, \, x>0\}$ with small data and study the scattering behaviour of solution.  To our knowledge there is no result on the well posedness of initial boundary value problem for Gross-Neveu model with general initial data in $L^2$. Our purpose is to prove the existence and the uniqueness in $C^1(\overline{\Omega})$ and in $L^2(\Omega)$ of global solution to (\ref{eq-dirac}-\ref{eq-dirac-bc}).

The first result is the following.

\begin{theorem}\label{thm-1}
Suppose that (H1) and (H2) hold. Let $(u_0,v_0)\in C^1([0,\infty))$ with  compact support in $[0,\infty)$ and satisfy the compatibility conditions as follow,
\begin{equation}\label{eq-compatcond1} u_0(0)=\lambda(0) v_0(0) \end{equation}
and
\begin{eqnarray}
&&(1-z_t(0))u_{0x}(0)+\lambda(0)(1+z_t(0)) v_{0x}(0)+i\lambda(0)\big(mu_0(0)-N_2(u_0(0),v_0(0))\big)  \nonumber \\
&&\qquad-i\big(mv_0(0)-N_1(u_0(0),v_0(0))\big)+\lambda_t(0)v_0(0)=0. \label{eq-compatcond2}
 \end{eqnarray}
Then (\ref{eq-dirac}-\ref{eq-dirac-bc}) has a unique global solution $(u,v)\in C^1(\overline{\Omega})$.
\end{theorem}

This result could be generalized to the following case.

\begin{theorem}\label{thm-2}
Suppose that (H1) and (H2) hold. Let $(u_0,v_0)\in H^1([0,\infty))$ satisfy the compatibility conditions as follows,
\begin{equation}\label{eq-compatcond3} u_0(0)=\lambda(0) v_0(0). \end{equation}  Then (\ref{eq-dirac}-\ref{eq-dirac-bc}) has a unique global solution $(u,v)\in H^1_{loc}(\Omega) \cap C(\overline{\Omega})$. Moreover,
\[ (u(\cdot,t),v(\cdot,t))\in H^1([z(t),\infty)) \]
for  $t\in [0,\infty)$.
\end{theorem}

With Theorem \ref{thm-1}, we can look for the global strong solution. Here the strong solution is defined as follows.

\begin{definition}
A pair of measurable functions $(u,v)$ is called a strong solution to (\ref{eq-dirac}-\ref{eq-dirac-bc}) if there exits a sequence of classical solutions $ (u^{(n)},v^{(n)})\in C^1(\overline{\Omega})$ to (\ref{eq-dirac}) such that
\[ u^{(n)}(z(t),t)=\lambda(t)v^{(n)}(z(t),t), \quad for  \, t>0,\]
and
\[  \lim_{n\to \infty}\int_{0}^N \big( |u^{(n)}(x,0)-u_0|^2+|v^{(n)}(x,0)-v_0|^2\big)dx=0,\]
\[ \lim_{n\to \infty}\iint_{K} \big( |u^{(n)}-u|^2+|v^{(n)}-v|^2\big)dxdt=0\]
for any compact set $K\subset \overline{\Omega}$ and for any $N>0$.
\end{definition}

\begin{theorem}\label{thm-3}
Suppose that (H1) and (H2) hold. For any $(u_0,v_0)\in L^2_{loc}([0,\infty))$, (\ref{eq-dirac}-\ref{eq-dirac-bc}) has a unique global strong solution $(u,v)\in L^2_{loc}(\Omega)$. Moreover, $|u||v|\in L^2_{loc}(\Omega)$, and $(u,v)$ solves  (\ref{eq-dirac}-\ref{eq-dirac-bc}) in the following sense,
\begin{eqnarray}
\mathop{\int\int}_{\Omega} \Big( i u(\phi_t+\phi_x)-mv\phi +N_1(u,v)\phi\Big) dxdt =-i \int^{\infty}_{0} u_0\phi(x,0)dx, \label{eq-weaksolu-1}\\
\mathop{\int\int}_{\Omega} \Big( i v(\psi_t-\psi_x)-mu\psi +N_2(u,v)\psi\Big) dxdt = -i \int^{\infty}_{0} v_0\psi(x,0)dx \label{eq-weaksolu-2}
\end{eqnarray}
for any $(\phi, \psi)\in C^1(\overline{\Omega})$ with bounded support in $\overline{\Omega}$ and $(\phi, \psi)(z(t),t)=0$ for $t\ge 0$.
\end{theorem}

Moreover, we have the following.

\begin{theorem}\label{thm-4}
Suppose that (H1) and (H2) hold. If $(u_0,v_0)\in L^2([0,\infty))$, then the strong solution $(u,v)$ given by Theorem \ref{thm-3} satisfies the following,
\[ (u,v)\in L^2(\Omega \cap (R^1\times [0,T])), \quad |u||v|\in L^2(\Omega \cap (R^1\times [0,T])) \] for any $T>0$.
Moreover, if $|\lambda(t)|^2(1-z_t(t))=(1+z_t(t))$ for $t\ge 0$, then
\[ \int_{z(t)}^{\infty} (|u(x,t)|^2+|v(x,t)|^2)dx=\int^{\infty}_0(|u_0(x)|^2+|v_0(x)|^2)dx \]
for almost every $t\in [0,\infty)$.
\end{theorem}

The remaining is organized as follows. First, in section \ref{section-existence-classical-solu}, to prove Theorem \ref{thm-1} and Theorem \ref{thm-2} for  (\ref{eq-dirac}-\ref{eq-dirac-bc}), we derive the equations (\ref{eq-dirac1}) and (\ref{eq-dirac2}) for $|u|^2$ and $|v|^2$ for local smooth solution $(u,v)$, and apply the characteristic method to the equations (\ref{eq-dirac1}) and (\ref{eq-dirac2}) to get the pointwise bounds on $|u|^2$ and $|v|^2$. Then it enables us to get the uniform $L^{\infty}$ bounds on $(u,v)$ in the domain $\Omega\cap \{0\le t <T \}$  for any $T>0$ and  extend the local solution globally.
 In section \ref{section-etimates-smooth-solu} we introduce  a Bony type functional $Q_0(t, \Delta)$ and a Glimm type functional $F_0(t,\Delta)=L(t,u,\Delta)+K_0L(t,v,\Delta)+C_0Q(t,\Delta)$ for smooth solution $(u,v)$ to get $L^2$ estimates of nonlinear term, $\int\int_{\Delta}|u(x,t)|^2|v(x,t)|^2dxdt$ on each characteristic triangle $\Delta$. Here different from the work in \cite{Zhang-zhao}, for the case that $\Delta\cap \partial\Omega\neq \emptyset$, by the assumption (H2) we choose a suitable constant $K_0>0$ so that the derivative of the weighted $L^2$ norm, $\frac{d}{dt}\big(L(t,u,\Delta)+K_0L(t,v,\Delta)\big)$ can control the possible increasing of the functional $Q_0(t,\Delta)$, and choose a suitable constant $C_0$ so that $F_0(t,\Delta)$ can control $\int\int_{\Delta}|u(x,t)|^2|v(x,t)|^2dxdt$,    while  for the case that $\Delta\subseteq \Omega$ same argument as in \cite{Zhang-zhao} can be carried out to get the control on $\int\int_{\Delta}|u(x,t)|^2|v(x,t)|^2dxdt$. In section \ref{section-difference-smooth-solu}, we consider the difference $(U,V)=(u-u^{\prime},v-v^{\prime})$ for two smooth solutions $(u,v)$ and $(u^{\prime},v^{\prime})$. We first write down the equations (\ref{eq-dirac3}) and (\ref{eq-dirac4}) for $(U,V)$, which contain $(U,V)$, $(u,v)$ and $(u^{\prime}, v^{\prime})$.  Then we introduce a Bony type functional $Q_1(t,\Delta)$ and a Glimm type functional $F_1(t,\Delta)$ for $|U|^2$, $|V|^2$, $|u|^2, |v|^2$ and $|u^{\prime}|^2$ and $|v^{\prime}|^2$,  and use it to prove the $L^2$ stability estimates in Proposition \ref{lemma-stability}. Here, as in section \ref{section-etimates-smooth-solu}, for the case that $\Delta\cap \partial\Omega\neq \emptyset$, by the assumption (H2) we choose a suitable constant $K>0$ so that the derivative of the weighted $L^2$ norm, $\frac{d}{dt}\big(L(t,U,\Delta)+K_1L(t,V,\Delta)\big)$ can control the possible increasing of the functional $Q_1(t,\Delta)$. In section \ref{section-convergence}, we first approximate the initial data (\ref{eq-dirac-initialv}) by  a sequence of smooth functions. Then,
  by the result on the global wellposedness for smooth solution in section \ref{section-existence-classical-solu}, we can have a sequence of global smooth solutions for smooth data for (\ref{eq-dirac}). With the help of the $L^2$ stability estimates in section \ref{section-difference-smooth-solu}, we show that the sequence of global smooth solutions converges to a strong solution in $L^2(\Delta)$ for any triangle $\Delta$. In section \ref{section-proof}, we complete the proof of Theorem \ref{thm-3} and Theorem \ref{thm-4}.

\section{Global classical solution}\label{section-existence-classical-solu}

For $T>0$, denote \[ \Omega(T)=\{ (x,t)\, | \, z(t)\le x<\infty, \, 0\le t<T \}.\]
Classical theory on semilinear hyperbolic systems \cite{alinhac} gives the following local existence result (see also \cite{mizohata}).

\begin{lemma} Suppose that the compatibility conditions (\ref{eq-compatcond1}) and (\ref{eq-compatcond2}) hold. For any $(u_0,v_0)\in C^1([0,\infty))$ with compact support in $[0,\infty)$,
there exists a $T_*>0$ such that (\ref{eq-dirac}-\ref{eq-dirac-bc}) has a unique solution  $(u,v)\in C^1(\overline{\Omega(T_*)})$.
\end{lemma}

Our aim in this section is to extend the solution $(u,v)$ globally to $\overline{\Omega}$. To this end, let  $(u_0,v_0)\in C^1([0,\infty))$ with compact support and let $(u,v)\in C^1(\Omega(T))$ be the solution to (\ref{eq-dirac}-\ref{eq-dirac-bc}) for $T\ge T_*$, taking $(u_0,v_0)$ as its initial data,  we have to establish the estimates on  $||(u,v)||_{L^{\infty}(\Omega(T))}$ in the next. Here we assume that the compatibility conditions (\ref{eq-compatcond1}) and (\ref{eq-compatcond2}) hold for $(u_0,v_0)$.

Multiplying the first equation of (\ref{eq-dirac}) by $\overline{u}$   and the second equation by $\overline{v}$ gives
\begin{equation}\label{eq-dirac1}
(|u|^2)_t +(|u|^2)_x=2m\Re (i\overline{u}v)+2\Re (i\overline{N_1}u),
\end{equation}
and
\begin{equation}\label{eq-dirac2}
 (|v|^2)_t-(|v|^2)_x=2m\Re (iu\overline{v})+2\Re (i\overline{N_2}v),
\end{equation}
 which, together with the structure of nonlinear terms, leads to
\begin{equation}\label{eq-dirac-conserv} (|u|^2+|v|^2)_t +(|u|^2-|v|^2)_x=0.\end{equation}

For the nonlinear terms in the righthand side of (\ref{eq-dirac1}) and (\ref{eq-dirac2}), we have the following by direct computation.
\begin{lemma}\label{lemma-nonlinearterm1}
Let $r_0(x,t)=m(|u(x,t)|^2+|v(x,t)|^2)+8|\beta||u(x,t)|^2|v(x,t)|^2$. Then
there hold the followings,
\[\big|2m\Re (i\overline{u}v)+2\Re (i\overline{N_1}u) \big|\le r_0(x,t)\]
and
\[\big|2m\Re (i\overline{v}u)+2\Re (i\overline{N_2}v)\big|\le r_0(x,t).\]
\end{lemma}
And we have the estimates on the $L^2$ norm of the solution as follows.
\begin{lemma}\label{lemma-conservation-ineq} Let $E_0=\int^{\infty}_0 (|u_0(x)|^2+|v_0(x)|^2)dx$. Then
for any $t\in [0,T)$, there holds the following,
\begin{equation}
\int^{\infty}_{z(t)}(|u(x,t)|^2+|v(x,t)|^2)dx\le  E_0.
\end{equation}
\end{lemma}
{\it Proof.} By (\ref{eq-dirac-bc}) and (\ref{eq-dirac-conserv}), and by assumption (H2), we have
\begin{eqnarray*}
&\,& \frac{d}{dt}\int^{\infty}_{z(t)}(|u(x,t)|^2+|v(x,t)|^2)dx  \\
&=& |u(z(t),t)|^2(1-z_t(t))  - |v(z(t),t)|^2(1+z_t(t)) \\
&=&|v(z(t),t)|^2[|\lambda(t)|^2(1-z_t(t))  - (1+z_t(t))]\le 0,
\end{eqnarray*}
which gives the desired inequality and completes the proof. $\Box$

 We consider the characteristic triangles for $(u,v)$ in $\Omega(T)$. For any $a, b\in R^1$ with $a<b$ and for any $t_0\ge 0$, we denote
\[\Delta(a,b,t_0)=\{ (x,t)\big|\, a-t_0+t<x<b+t_0-t, \, t_0<t<\frac{b-a}{2}+t_0\}, \]
 see Figure \ref{fig-domain1},
\begin{figure}[ht]
\begin{center}
\unitlength=10mm
\begin{picture}(10,3.5)
\thicklines
\put(0,0){\line(1,0){10}}
\put(1,0){\line(5,4){4}}
\put(9,0){\line(-5,4){4}}
\put(0.5,-0.5){$(a,t_0)$}
\put(8.5,-0.5){$(b,t_0)$}\put(10.3,0){$t=t_0$}
\put(4.2,3.5){$(\frac{a+b}{2},\frac{b-a}{2}+t_0)$}
\put(4,1.5){$\Delta(a,b,t_0)$}
\end{picture}
\caption{Domain $\Delta(a,b,t_0)$}\label{fig-domain1}
\end{center}
\end{figure}
and, denote
\[ \Gamma_u(x_0,t_0;t_1)=\{(x,t)\, \big|\, x=x_0-t_0+t, \, t_1\le t \le t_0\} \]
and
\[ \Gamma_v(x_0,t_0;t_1)=\{(x,t)\, \big|\, x=x_0+t_0-t, \, t_1\le t\le t_0\} \]
for $t_1\le t_0$, see Figure \ref{fig-charact}.
\begin{figure}[ht]
\begin{center}
\unitlength=10mm
\begin{picture}(10,3.8)(0,-0.3)
\thicklines
\put(0,0){\line(1,0){10}}
\put(1,0){\line(5,4){4}}
\put(9,0){\line(-5,4){4}}
\put(0,-0.5){$(x_0-t_0+t_1,t_1)$}
\put(8,-0.5){$(x_0+t_0-t_1,t_1)$}
\put(10.3,0){$t=t_1$}
\put(4.5,3.5){$(x_0,t_0)$}
\put(3.3,1.5){$\Gamma_u(x_0,t_0;t_1)$}
\put(7.5,1.5){$\Gamma_v(x_0,t_0;t_1)$}
\end{picture}
\caption{Characteristic lines $\Gamma_u$ and $\Gamma_v$}\label{fig-charact}
\end{center}
\end{figure}
 It is obvious that $\Gamma_u(x_0,t_0;t_1)$ is a characteristic line for the first equation of $u$ in (\ref{eq-dirac}) while $\Gamma_v(x_0,t_0;t_1)$ is a characteristic line for the second equation of $v$ in (\ref{eq-dirac}).

 Along these characteristic lines in $\Omega(T)$, we have the following estimates.

\begin{lemma}\label{lemma-v-charact}
If $\Gamma_v(x_0,t_0; t_1)\subseteq\overline{\Omega(T)}$, then
\[ \int_{t_1}^{t_0} |u(x_0+t_0-s,s)|^2 ds \le  E_0.\]
Here $E_0=\int^{\infty}_0(|u_0|^2+|v_0|^2)dx$.
\end{lemma}
{\it Proof.} Denote
\[ \omega(x_0,t_0)=\{ (x,t)| z(t)\le x\le x_0+t_0-t, \, 0\le t \le t_0\}.\]
Then taking the integration of (\ref{eq-dirac-conserv}) over $\omega(x_0,t_0)$ gives the following,
\begin{eqnarray*}
&\, & \int^{x_0+t_0}_0(|u_0(x)|^2+|v_0(x)|^2)dx \\
&=&2\int_0^{t_0}|u(x_0+t_0-s,s)|^2 ds +\int_{z(t_0)}^{x_0} (|u(x,t_0)|^2+|v(x,t_0)|^2 )dx \\
&\,& +\int^{t_0}_0 \big\{ (-1+z_t(s))|u(z(s),s)|^2+ (1+z_t(s))|v(z(s),s)|^2\big\} ds \\ &\ge& \int_0^{t_0}|u(x_0+t_0-s,s)|^2 ds,
\end{eqnarray*}
where we use the boundary condition (\ref{eq-dirac-bc}) and assumption (H2) to get the last inequality. This implies the result and the proof is complete. $\Box$

\begin{lemma}\label{lemma-u-charact}
If $\Gamma_u(x_0,t_0;t_1)\subseteq \overline{\Omega(T)}$, then
\[  \int^{t_0}_{t_1} |v(x_0-t_0+s,s)|^2 ds\le E_0.\]
Here $E_0=\int^{\infty}_0(|u_0|^2+|v_0|^2)dx$.
\end{lemma}
{\it Proof.} Since $\Gamma_u(x_0,t_0;t_1)\subseteq \overline{\Omega(T)}$, then the domain \[\Delta(x_0-t_0+t_1, x_0+t_0-t_1,t_1)\subseteq\overline{\Omega(T)}.\]

Taking the integration of (\ref{eq-dirac-conserv})  over $\Delta(x_0-t_0+t_1, x_0+t_0-t_1,t_1)$, we have
\begin{eqnarray*}
\int^{t_0}_{t_1}|v(x_0-t_0+s,s)|^2ds &\le& \int^{x_0+t_0-t_1}_{x_0-t_0+t_1}(|u(x,t_1)|^2+|v(x,t_1)|^2)dx \\
&\le& \int_{z(t_1)}^{\infty} (|u(x,t_1)|^2+|v(x,t_1)|^2)dx \\
&\le&  E_0,
\end{eqnarray*}
where the last inequality is given by Lemma \ref{lemma-conservation-ineq}.
 The proof is complete. $\Box$

Using the above estimates on  along the characteristic lines, we can get the following pointwise estimates on $v$ at first.

\begin{lemma}\label{lemma-pointwise-v}
For $(x,t)\in \overline{\Omega(T)}$,
\[ |v(x,t)|^2 \le (|v_0(x+t)|^2+mE_0)\exp(mt+8|\beta| E_0).\]
Here $ E_0=\int_0^{\infty}(|u_0|^2+|v_0|^2)dx$.
\end{lemma}
{\it Proof.} Assumption (H1) implies that \[\Gamma_v(x,t; 0)\subset \overline{\Omega(T)}\] for any $(x,t)\in \overline{\Omega(T)}$.

Then, by Lemma \ref{lemma-nonlinearterm1}, along $\Gamma_v(x,t; 0)$ we use the  equation  (\ref{eq-dirac2}) to derive that
\[ \frac{d}{ds}|v(x+t-s,s)|^2 \le m|u(x+t-s,s)|^2+\big(m+8|\beta||u(x+t-s,s)|^2\big)|v(x+t-s,s)|^2.\]
Therefore
\begin{eqnarray*} \frac{d}{ds}\big(|v(x+t-s,s)|^2e_1(x,t,s)\big) &\le& m|u(x+t-s,s)|^2e_1(x,t,s) \\ &\le& m|u(x+t-s,s)|^2,
\end{eqnarray*}
where
\[ e_1(x,t,s)=\exp\big(-ms-8|\beta|\int_0^s|u(x+t-\tau,\tau)|^2d\tau\big).\]
Taking the integration of the above from $s=0$ to $t$, we can prove the desired result by Lemma \ref{lemma-v-charact}. The proof is complete. $\Box$

To get the pointwise estimates on $u$, we look for the intersection point of the boundary $\Gamma_B$ and the characteristic line $\{(x,t)|x-t=b\}$ for $b\le 0$.

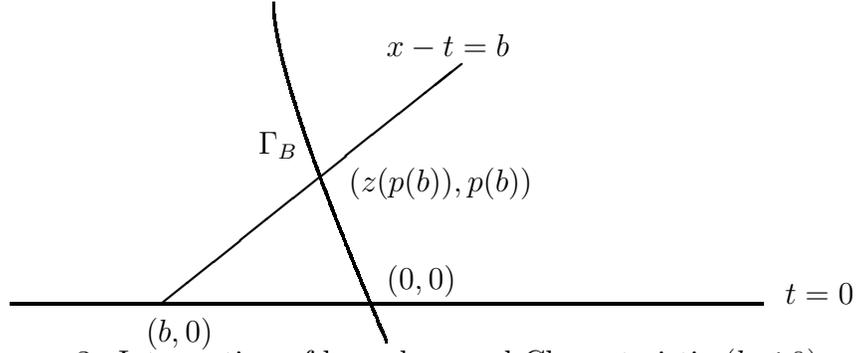
\begin{figure}[ht]
\begin{center}
\unitlength=10mm
\begin{picture}(10,3.5)
\thicklines
\put(0,0){\line(1,0){10}}
\put(2,0){\line(5,4){4}}
\put(5,3.3){$x-t=b$}
\put(1.8,-0.5){$(b,0)$}
\put(5,0.2){$(0,0)$}
\put(10.3,0){$t=0$}
\qbezier(5,-0.5)(3.5,3)(3.5,4)
\put(4.5,1.5){$(z(p(b)),p(b))$}
\put(3.3,2){$\Gamma_B$}
\end{picture}
\caption{Intersection of boundary and Characteristic $(b\le 0)$}\label{fig-intersect-bdry-charac}
\end{center}
\end{figure}
\begin{lemma}\label{lemma-geometry-bdry}
For any $b\le 0$, the equation $z(t)-t=b$ has a unique solution $t=p(b)$, where $p\in C^1 (-\infty,0]$ and $p^{\prime}(s)<0$ for $s\le 0$.
\end{lemma}
{\it Proof.} From assumption (H1) it follows that \[z_t(t)-1<0\]
 for $t>0$, which implies that the function $z(t)-t$ has a global inverse $p\in C^1(-\infty,0]$. Moreover, \[p^{\prime}(s)=\frac{1}{1-z_t(p(s))}<0.\] Therefore the proof is complete.$\Box$

Now we can have the following pointwise estimates on $u$.

\begin{lemma}\label{lemma-pointwise-u}
If $(x,t)\in \overline{\Omega(T)}$ with $x-t\ge 0$, then
\[ |u(x,t)|^2 \le (|u_0(x-t)|^2+mE_0)\exp(mt+8|\beta|E_0).\]
If $(x,t)\in \overline{\Omega(T)}$ with $x-t< 0$, then
\[ |u(x,t)|^2  \le (|\lambda(t)|^2+1)(|v_0(2p(x-t)+x-t)|^2 +mE_0) \exp(2mt+16|\beta|E_0).\]
Here $E_0=\int^{\infty}_0(|u_0|^2+|v_0|^2)dx$.
\end{lemma}
{\it Proof.} For $(x,t)\in \overline{\Omega(T)}$ with $x-t\ge 0$, the assumption (H1) implies that
\[\Gamma_u(x,t;0)\subseteq \overline{\Omega(T)}.\] Then, by (\ref{eq-dirac1}) and by Lemma \ref{lemma-nonlinearterm1},
  we have
\begin{eqnarray}
\frac{d}{ds}\big( |u(x-t+s,s)|^2 e_2(x,t,s)\big)  &\le&
m|v(x-t+s,s)|^2e_2(x,t,s) \nonumber\\ &\le& m|v(x-t+s,s)|^2, \label{eq-charact-u}
\end{eqnarray}
where \[ e_2(x,t,s)=\exp(-ms-8|\beta|\int^s_0 |v(x-t+\tau,\tau)|^2 d\tau).\]
Taking the integration of (\ref{eq-charact-u}) from $0$ to $t$ and using Lemma \ref{lemma-u-charact}, we get
\begin{eqnarray*}
|u(x,t)|^2
\le (|u_0(x-t)|^2 + mE_0)\exp(mt+8|\beta|E_0).
\end{eqnarray*}

For $(x,t)\in \overline{\Omega(T)}$ with $x-t< 0$, Lemma \ref{lemma-geometry-bdry} implies that the characteristic line $\Gamma_u(x,t;0)$ and the boundary intersect only at the  point $(z\big(p(x-t)\big),p(x-t))$.

Then, by (\ref{eq-dirac1}) and by Lemma \ref{lemma-nonlinearterm1}, along the characteristic line $\Gamma_u(x,t; p(x-t))$ we have
\begin{eqnarray*}
\frac{d}{ds}\big( |u(x-t+s,s)|^2 e_3(x,t,s) \big) & \le&
m|v(x-t+s,s)|^2e_3(x,t,s) \\ &\le& m|v(x-t+s,s)|^2 .
\end{eqnarray*}
where
\[ e_3(x,t,s)=\exp\big(-m(s-p(x-t))-8|\beta|\int^s_{p(x-t)} |v(x-t+\tau,\tau)|^2 d\tau\big).\]

Taking the integration of the above from $p(x-t)$ to $t$,  we use Lemma \ref{lemma-u-charact} and Lemma \ref{lemma-pointwise-v} to get the following,
\begin{eqnarray*}
|u(x,t)|^2 &\le& (|u(p(x-t)+x-t, p(x-t))|^2 +mE_0) \exp(mt+8|\beta|E_0) \\
&\le& (|\lambda(t)|^2|v(p(x-t)+x-t, p(x-t))|^2 +mE_0) \exp(mt+8|\beta|E_0) \\
&\le& (|\lambda(t)|^2+1)(|v_0(2p(x-t)+x-t)|^2 +mE_0) \exp(2mt+16|\beta|E_0).
\end{eqnarray*}
The proof is complete. $\Box$

Now using the pointwise estimates on $u$ and $v$, we can prove Theorem \ref{thm-1}.

{\bf Proof of Theorem \ref{thm-1}.} For $(u_0,v_0)\in H^1([0,\infty))\subset L^{\infty}([0,\infty))$, Lemma \ref{lemma-pointwise-v} and Lemma \ref{lemma-pointwise-u} lead to
\[ |u(x,t)|^2+|v(x,t)|^2\le (|\lambda(t)|^2+2)\big(||(u_0,v_0)||_{L^{\infty}}+2mE_0\big)\exp(2mt+16|\beta|E_0)\]
for $x\ge z(t)$ and $0\le t<T$.

Then by the standard theory on semilinear hyperbolic equations (see \cite{alinhac} for instance), we can extend the solution $(u,v)$ across the time $t=T$.

Therefore, repeating the same argument for any time, we can extend the solution globally to $\overline{\Omega}$. The proof is complete. $\Box$

Furthermore Theorem \ref{thm-2} follows from Theorem \ref{thm-1}.

{\bf Proof of Theorem \ref{thm-2}.}  Let $(\phi_0,\psi_0)\in C^{\infty}_c(R^1)$ be a pair of functions such that $\phi_0(x)=u_0(0)$ and $\psi_0(x)=v_0(0)$ for $x$ belonging to a neighbourhood of zero. Then we choose a sequence of functions $(\phi_1^k,\psi_1^k)\in C^{\infty}_c(0,\infty)$ such that
 $(u_0^{(k)},v_0^{(k)}):=(\phi_0+\phi_1^k,\psi_0+\psi_1^k)$ is convergent to $(u_0,v_0)$ in $H^1(0,\infty)$ as $k$ tends to $\infty$.

 It is obvious that  $(u_0^{(k)},v_0^{(k)})$ satisfies the compatibility conditions as (\ref{eq-compatcond1}) and (\ref{eq-compatcond2}). Therefore, by Theorem \ref{thm-1}, the equations (\ref{eq-dirac})  has a global smooth solution $(u^{(k)},v^{(k)})$  with the initial data $(u_0^{(k)},v_0^{(k)})$ for $k\ge 1$.

 Moreover, by Lemma \ref{lemma-pointwise-v} and Lemma \ref{lemma-pointwise-u}, we have
 \[ \sup_{k\ge 1} ||(u^{(k)},v^{(k)})||_{L^{\infty}(\Omega(T))}<\infty\]
 for any $T>0$, which enables us to show as in \cite{alinhac} and \cite{mizohata} that the sequence $(u^{(k)},v^{(k)})$ is convergent in $H^1(\Omega(T))$ to a solution $(u,v)$ of (\ref{eq-dirac})-(\ref{eq-dirac-bc}) as $k$ tends to $\infty$ for any $T>0$.

  The uniqueness can be proved by the the  energy inequality for the difference of solutions in $L^{\infty}(\Omega(T))\cap H^1(\Omega(T))$ as in \cite{alinhac} and \cite{mizohata}.  The proof is complete.$\Box$

\section{ Estimates on the classical solution}\label{section-etimates-smooth-solu}

Consider the case that $(u_0,v_0)\in C^1([0,\infty))$, and let $(u,v)\in C^1 (\overline{\Omega})$ be the global solution to (\ref{eq-dirac}) with boundary condition (\ref{eq-dirac-bc}). Here we assume that the compatibility condition (\ref{eq-compatcond1}) and (\ref{eq-compatcond2}) hold.  Our aim in this section is to establish the local estimates on  $(u,v)$.

To this end, set $\Delta=\Delta(a,b;t_0)$ for simplification  and assume that  $\Delta\cap \Omega \neq \emptyset$ in this section.

Let
\[ x_0^{\prime}=\frac{b+a}{2}, \, \, \, t_0^{\prime}=\frac{b-a}{2}+t_0.\]
Then $\Gamma_u(x_0^{\prime},t_0^{\prime}; t_0)$ and $\Gamma_v(x_0^{\prime},t_0^{\prime}; t_0)$ are the left and right edges of $\Delta$. By (H1), $\Gamma_B$ and $\Gamma_u(x_0^{\prime},t_0^{\prime}; t_0)\cup\Gamma_v(x_0^{\prime},t_0^{\prime}; t_0)$ intersect at one point at the most.

We introduce a time interval as follows.
Denote \[ I_{\Delta}=\{t\big|\, t_0\le t\le \frac{b-a}{2}+t_0,\, z(t)\le b+t_0-t\}.\]

By Lemma \ref{lemma-geometry-bdry}, we have the following.

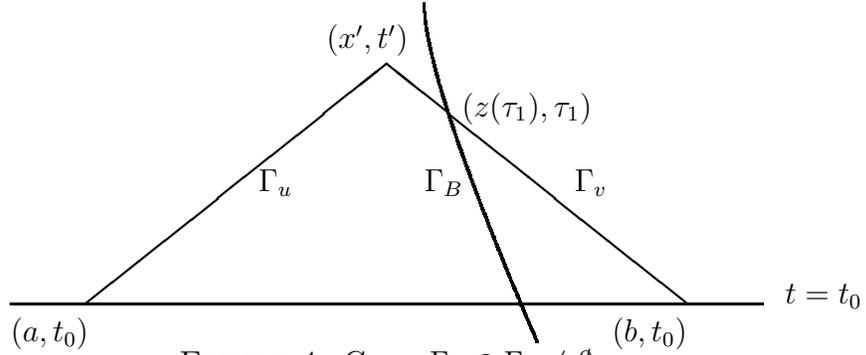
\begin{figure}[ht]
\begin{center}
\unitlength=10mm
\begin{picture}(10,3.5)
\thicklines
\put(0,0){\line(1,0){10}}
\put(1,0){\line(5,4){4}}
\put(9,0){\line(-5,4){4}}
\put(0,-0.5){$(a,t_0)$}\put(5.5,1.5){$\Gamma_B$}
\put(8,-0.5){$(b,t_0)$}\put(10.3,0){$t=t_0$}
\put(3.3,1.5){$\Gamma_u$}\put(7.5,1.5){$\Gamma_v$}
\qbezier(7,-0.5)(5.5,3)(5.5,4)
\put(6,2.5){$(z(\tau_1),\tau_1)$}
\put(4.2,3.4){$(x^{\prime},t^{\prime})$}
\end{picture}
\caption{Case: $\Gamma_B\cap\Gamma_v\neq\emptyset$}\label{fig-case1}
\end{center}
\end{figure}

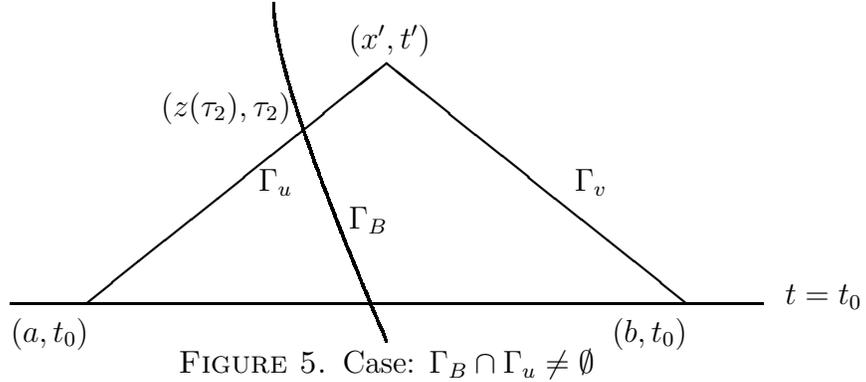
\begin{figure}[ht]
\begin{center}
\unitlength=10mm
\begin{picture}(10,3.5)
\thicklines
\put(0,0){\line(1,0){10}}
\put(1,0){\line(5,4){4}}
\put(9,0){\line(-5,4){4}}
\put(0,-0.5){$(a,t_0)$}
\put(8,-0.5){$(b,t_0)$}\put(10.3,0){$t=t_0$}
\qbezier(5,-0.5)(3.5,3)(3.5,4)
\put(2,2.5){$(z(\tau_2),\tau_2)$}
\put(4.5,3.4){$(x^{\prime},t^{\prime})$}
\put(3.3,1.5){$\Gamma_u$}\put(7.5,1.5){$\Gamma_v$}\put(4.5,1){$\Gamma_B$}
\end{picture}
\caption{Case: $\Gamma_B\cap\Gamma_u\neq\emptyset$}\label{fig-case2}
\end{center}
\end{figure}

\begin{lemma}\label{lemma-timeinterval} There hold the following statements. (1) If
$\Gamma_B\cap\Gamma_v(x_0^{\prime},t_0^{\prime}; t_0)=\{(z(\tau_1),\tau_1)\}$ for some $\tau_1\ge t_0$ (Figure \ref{fig-case1}), then $I_{\Delta}=[t_0,\tau_1]$ and
\[ \{x\big|\, (x,t)\in \overline{\Delta\cap\Omega}\}=[z(t),b+t_0-t].\]
(2) If
$\Gamma_B\cap\Gamma_u(x_0^{\prime},t_0^{\prime}; t_0)=\{(z(\tau_2),\tau_2)\}$ for some $\tau_2\ge t_0$  (Figure \ref{fig-case2}),  then $I_{\Delta}=[t_0,\frac{b-a}{2}+t_0]$ and
\[ \{x\big|\, (x,t)\in \overline{\Delta\cap\Omega}\}=[z(t),b+t_0-t] \, \, \, \mbox{for} \,  t\in [t_0, \tau_2],\]
and
\[ \{x\big|\, (x,t)\in \overline{\Delta\cap\Omega}\}=[a-t_0+t,b+t_0-t] \, \, \, \mbox{for} \,  t\in [ \tau_2, \frac{b-a}{2}+t_0].\]
(3) If $\Delta\subset \Omega$, then $I_{\Delta}=[t_0,\frac{b-a}{2}+t_0]$
and
\[ \{x\big|\, (x,t)\in \overline{\Delta\cap\Omega}\}=[a-t_0+t,b+t_0-t].\]
\end{lemma}

Now we can define the functionals for $(u,v)$ on $\Delta\cap\Omega$ as follow.

\begin{definition}\label{def-L}
For $t\in I_{\Delta}$, and for any $w\in C^1(\overline{\Omega})$, define,
\begin{equation} L(t,w, \Delta)=
\int_{z_a(t)}^{b-t+t_0} |w(x,t)|^2 dx,
\end{equation}
where \[z_a(t)= \max\{ a+t-t_0,z(t)\}.\]
\end{definition}

\begin{definition}\label{def-L0D0Q0}
For $t\in I_{\Delta}$,
and for the solution $(u,v)$, define
\[
L_0(t,\Delta)=L(t,u,\Delta)+L(t,v,\Delta)
\]
and
\[
D_0(t,\Delta)=\int_{z_a(t)}^{b-t+t_0} |u(x,t)|^2|v(x,t)|^2 dx,
\]
\[ Q_0(t,\Delta)=
\int\int_{z_a(t)<x<y<b-t+t_0} |u(x,t)|^2|v(y,t)|^2 dxdy,
\]
where \[z_a(t)= \max\{ a+t-t_0,z(t)\}.\]
\end{definition}

Then we have the following estimates on the $L^2$ norm.

\begin{lemma}\label{lemma-local-charge}
For $t\in I_{\Delta}$, there holds the following,
\[ L_0(t,\Delta)\le L_0(t_0,\Delta).\]
\end{lemma}

{\it Proof.} It suffices to prove lemma for three cases according to Lemma \ref{lemma-timeinterval}.

Case 1: The right edge of  $\Delta$ and $\Gamma_B$ intersect at some point $(z(\tau_1),\tau_1)$, see Figure \ref{fig-case1}. In this case $I_{\Delta}=[t_0,\tau_1]$.

 Then for $t\in [t_0,\tau_1]$, $z_a(t)=z(t)$. Moreover, by (\ref{eq-dirac-bc}) and (\ref{eq-dirac-conserv}), and by assumption (H2), we have
\begin{eqnarray*}
&\,& \frac{d}{dt}\int^{b-t+t_0}_{z(t)}(|u(x,t)|^2+|v(x,t)|^2)dx  \\ &=& |u(z(t),t)|^2(1-z_t(t))  - |v(z(t),t)|^2(1+z_t(t)) \\
&\,& -(|u(b-t+t_0,t)|^2+|v(b-t+t_0,t)|^2) \\
&\le&|v(z(t),t)|^2[|\lambda(t)|^2(1-z_t(t))  - (1+z_t(t))]\le 0.
\end{eqnarray*}
This leads to the desired result.

Case 2: The left edge of $\Delta$ and $\Gamma_B$ intersect at some point $(z(\tau_2),\tau_2)$, see Figure \ref{fig-case2}. Then $I_{\Delta}=[t_0,t_0+\frac{b-a}{2}]$.

For $t\in [t_0,\tau_2]$, $z_a(t)=z(t)$, and in the same way as in the proof of Case 1, we can get
\[  L_0(t,\Delta)\le L_0(t_0,\Delta).\]
For $t\in [\tau_2,\frac{b-a}{2}+t_0]$, $z_a(t)=a-t_0+t$, then we can use the result for Case 2 to deduce that
\[  L_0(t,\Delta)\le L_0(\tau_2,\Delta).\]

Case 3: $\Delta$ lies in the interior of $\Omega$. The proof can be carried out in the same way as in Case 1.

Therefore the proof is complete. $\Box$

For any $T>0$, we recall the notation
 \[ \Omega(T)=\{ (x,t)\, | \, z(t)\le x<\infty, \, 0\le t<T \},\]
 and have the control on the potential $Q_0$ for the case that $\Delta \subset \Omega(T)$ as follows.

\begin{lemma}\label{lemma-int-Q}Suppose that $\Delta \subset \Omega(T)$ for $T>0$. Then there exists  constants $\delta_0>0$ such that for the initial data satisfying $L_0(t_0,\Delta)\le \delta_0$ there holds the following
\begin{eqnarray}
\frac{d Q_0(t,\Delta)}{dt }+D_0(t,\Delta)\le 2m(L_0(t_0,\Delta))^2
\end{eqnarray}
for $t\in (t_0, \frac{b-a}{2}+t_0)$.
Therefore,
\begin{eqnarray}
 Q_0(t,\Delta) + \int\limits_{t_0}^t D_0(\tau,\Delta) d\tau &\le& 2m(L_0(t_0,\Delta))^2 (t-t_0) + Q_0(t_0,\Delta)  \,  \,  \,  \nonumber \\ &\le& 2m(L_0(t_0,\Delta))^2 (t-t_0)+ (L_0(t_0,\Delta))^2 \,  \,  \,
\end{eqnarray}
for $t\in [t_0, \frac{b-a}{2}+t_0]$. Here $\delta_0$ is independent of $T$.
\end{lemma}

The proof of Lemma \ref{lemma-int-Q} has been given in \cite{Zhang-zhao} and is similar to the proof of Lemma \ref{lemma-bdry-Q} in the next.

To get the control on the potential $Q_0$ near the boundary, we introduce a new functional as follows.
\begin{definition}\label{def-F0} For constants $K_0>0$ and $C_0>0$ and for $t\in I_{\Delta}$, define
\[  F_0(t,\Delta)=L(t,u,\Delta)+K_0L(t,v,\Delta) +C_0Q_0(t,\Delta). \]
\end{definition}

For any $T>0$, we have the control on $F_0$ near the boundary as follows.
\begin{lemma}\label{lemma-bdry-Q}
Suppose that $\Delta\subset R^1\times [0,T]$ and $\Delta\cap \Gamma_B\neq \emptyset$ for $T>0$. Then there exist constants  $\delta_0>0$, $K_0>0$ and $C_0>0$ such that for $L_0(t_0,\Delta)\le \delta_0$ there hold the following,
\begin{equation}\label{eq-ineqF0}
\frac{d}{dt} F_0(t,\Delta)\le -D_0(t,\Delta)-|v(z_a(t),t)|^2 +O(1)\delta_0,
\end{equation}
for $t\in I_{\Delta}$ with $z(t)\neq a+t-t_0$.
 Here the constants $\delta_0,K_0$ and $C_0$ depend only on $T$; and the bound of $O(1)$ depends only on $T$.
\end{lemma}

{\it Proof.}  For simplification, we denote $L_0(t, \Delta)$,$D_0(t,\Delta)$, $F_0(t,\Delta)$ and $Q_0(t,\Delta)$  by $L_0(t)$,$D_0(t)$, $F_0(t)$ and $Q_0(t)$. Now it suffices to prove the lemma for two cases.

 Case 1: The boundary $\Gamma_B$ and the right edge $\Gamma_v$ of $\Delta$  intersect at the point $(z(\tau_1),\tau_1)$ for some $\tau_1\in [t_0, t_0+\frac{b-a}{2}]$, see Figure \ref{fig-case1}.

 Then  $I_{\Delta}= [t_0,\tau_1]$, $z_a(t)=z(t)$. For $t\in [t_0,\tau_1]$, by Lemma \ref{lemma-nonlinearterm1}, we use (\ref{eq-dirac1}), (\ref{eq-dirac2})  to get
\begin{eqnarray*}
\frac{d}{dt} L(t,u,\Delta)&\le & \int^{b-t+t_0}_{z(t)}\big(-(|u(x,t)|^2)_x+r_0(x,t)\big) dx \\
&\,& -|u(b-t+t_0,t)|^2-z_t(t)|u(z(t),t)|^2 \\
&\le& (1-z_t(t))|u(z(t),t)|^2+\int^{b-t+t_0}_{z(t)}r_0(x,t)dx,
\end{eqnarray*}
and
\begin{eqnarray*}
\frac{d}{dt} L(t,v,\Delta)&\le & \int^{b-t+t_0}_{z(t)}\big((|v(x,t)|^2)_x+r_0(x,t)\big)dx \\
&\,& -|v(b-t+t_0,t)|^2-z_t(t)|v(z(t),t)|^2 \\
&\le& -(1+z_t(t))|v(z(t),t)|^2+\int^{b-t+t_0}_{z(t)}r_0(x,t)dx,
\end{eqnarray*}
which lead to
\begin{eqnarray*}
\frac{d}{dt}\big( L(t,u,\Delta)+K_0L(t,v,\Delta)\big) &\le & \big((1-z_t(t))|\lambda(t)|^2-K_0(1+z_t(t))\big)|v(z(t),t)|^2 \\ &\,& +(1+K_0)\int^{b-t+t_0}_{z(t)}r_0(x,t)dx
\\ &\le & -2|v(z(t),t)|^2+(1+K_0)\int^{b-t+t_0}_{z(t)}r_0(x,t)dx \\ &\le & -2|v(z(t),t)|^2+ O(1)(L_0(t)+D_0(t)),
\end{eqnarray*}
where we choose $K_0>1$ large enough so that
\begin{equation}\label{eq-K} (1-z_t(t))|\lambda(t)|^2-K_0(1+z_t(t))<-2 \end{equation} for $t\in [0,T]$.

On the other hand, by Lemma \ref{lemma-nonlinearterm1}, we use (\ref{eq-dirac1}), (\ref{eq-dirac2})  again to get the following for $Q_0$,

\begin{eqnarray*}
\frac{d}{dt}Q_0(t)&=&\int\int_{z(t)<x<y<b-t+t_0} \big( (|u(x,t)|^2)_t|v(y,t)|^2+|u(x,t)|^2(|v(y,t)|^2)_t\big) dxdy \\
& \, & +\big(\frac{d}{dt} \int\int_{z_a(t)<x<y<b-t+t_0} |u(x,s)|^2|v(y,s)|^2 dxdy \big)\big|_{s=t} \\
&\le & \int_{z(t)}^{b-t+t_0} (|u(z(t),t)|^2-|u(y,t)|^2)|v(y,t)|^2 dy   \\
&\,& +\int_{z(t)}^{b-t+t_0} |u(x,t)|^2(|v(b-t+t_0,t)|^2-|v(x,t)|^2)dx \\
&\,& +\int_{z(t)}^{b-t+t_0}r_0(x,t)dx\int_{z(t)}^{b-t+t_0}(|u(y,t)|^2+|v(y,t)|^2)dy \\
&\,& -z_t(t)\int^{b-t+t_0}_{z(t)} |u(z(t),t)|^2|v(y,t)|^2 dy \\
&\,& -\int_{z(t)}^{b-t+t_0} |u(x,t)|^2|v(b-t+t_0,t)|^2dx \\
 &\le& (-2+O(1)L_0(t))\int_{z(t)}^{b-t+t_0} |u(x,t)|^2|v(x,t)|^2 dx \\
 &\,& +(1-z_t(t))|u(z(t),t)|^2 \int_{z(t)}^{b-t+t_0}|v(y,t)|^2 dy+ O(1) (L_0(t))^2.
\end{eqnarray*}

Therefore,
\begin{eqnarray*}
\frac{d}{dt} F_0(t) & \le & -D_0(t)\big( (2-O(1)L_0(t))C_0-O(1)\big) \\
&\,& -|v(z(t),t)|^2\big(2-C_0(1-z_t(t))|\lambda(t)|^2 L_0(t)\big) \\
&\,& +L_0(t)(O(1)+O(1)C_0L_0(t)) \\
&\le& -D_0(t)-|v(z(t),t)|^2+O(1)L_0(t) \\ &\le& -D_0(t)-|v(z(t),t)|^2+O(1)L_0(t_0) ,
\end{eqnarray*}
where we choose constant $C_0>0$ and $\delta_0$ such that $ L_0(t_0)\le \delta_0$ and
\[ (2-O(1)\delta_0)C_0-O(1)\ge 1,\]
\[ 2-C_0(1-z_t(t))|\lambda(t)|^2\delta_0\ge 1\]
for $t\in [0,T]$. Then (\ref{eq-ineqF0}) is proved for this case.

Case 2: The boundary $\Gamma_B$ and the right edge $\Gamma_v$ of $\Delta$  intersect at the point $(z(\tau_2),\tau_2)$ for some $\tau_2\in [t_0, t_0+\frac{b-a}{2}]$, see Figure \ref{fig-case2}.
The proof of (\ref{eq-ineqF0}) can be carried out in the same way as in Case 1 for $t\neq\tau_2$.
Thus the proof is complete. $\Box$

\section{Estimates on the difference between the classical solutions}\label{section-difference-smooth-solu}

Let $(u,v)\in C^1(\overline{\Omega})$ and $(u^{\prime},v^{\prime})\in C^1(\overline{\Omega})$ be two classical solutions to (\ref{eq-dirac}) with (\ref{eq-dirac-bc}). We consider the difference between these two solutions  and denote \[(U,V)=(u-u^{\prime},v-v^{\prime}).\] Then,
\begin{eqnarray*}
U_t+U_x=imV-i(N_1(u,v)-N_1(u^{\prime},v^{\prime})), \\
V_t-V_x=imU-i(N_2(u,v)-N_2(u^{\prime},v^{\prime})),
\end{eqnarray*}
which lead to
\begin{eqnarray}
(|U|^2)_t+(|U|^2)_x=\Re 2\{ imV\overline{U}-i(N_1(u,v)-N_1(u^{\prime},v^{\prime}))\overline{U}\}, \label{eq-dirac3} \\
(|V|^2)_t-(|V|^2)_x=\Re 2\{imU\overline{V}-i(N_2(u,v)-N_2(u^{\prime},v^{\prime}))\overline{V}\}, \label{eq-dirac4}
\end{eqnarray}

For the nonlinear terms in the righthandsides of (\ref{eq-dirac3}) and (\ref{eq-dirac4}), we have following by direct computations.
\begin{lemma}\label{lemma-nonlinearterm2}
There exists a $c_*>0$ such that
\[ |\Re 2\{ imV\overline{U}-i(N_1(u,v)-N_1(u^{\prime},v^{\prime}))\overline{U}\}|\le r_1(x,t)\]
and
\[ |\Re 2\{imU\overline{V}-i(N_2(u,v)-N_2(u^{\prime},v^{\prime}))\overline{V}\}|\le r_1(x,t),\]
where
\[ r_1(x,t)=m(|U(x,t)|^2+|V(x,t)|^2)+ c_* r_2(x,x,t),\]
\[ r_2(x,y,t)=|U(x,t)|^2\big( |v(y,t)|^2+|v^{\prime}(y,t)|^2\big) +\big( |u(x,t)|^2+|u^{\prime}(x,t)|^2\big)|V(y,t)|^2.\]
\end{lemma}

To get the control on $(U,V)$ via (\ref{eq-dirac3}) and (\ref{eq-dirac4}), we introduce following functionals on $\Delta\cap\Omega$ for $(U,V)$ as in \cite{Zhang-zhao}. Here it is assume that $\Delta\cap \Omega\neq\emptyset$.
\begin{definition}
For $\Delta=\Delta(a,b,t_0)$ and $K_1>0$, $C_1>0$, define
\[ L_1(t,\Delta)=L(t,U,\Delta)+K_1 L(t,V,\Delta),\]
\[ D_1(t,\Delta)=\int_{z_a(t)}^{b-t+t_0} r_2(x,x,t)dx,\]
\[ Q_1(t,\Delta)=\int\int_{z_a(t)<x<y<b-t+t_0} r_2(x,y,t)dxdy\]
and
\[ F_1(t,\Delta)=L_1(t,\Delta)+C_1Q_1(t,\Delta)\]
for $t\in I_{\Delta}$.
Here  $\Delta\cap \Omega\neq\emptyset$ with \[z_a(t)= \max\{ a+t-t_0,z(t)\}; \]  $L(t,U,\Delta)$,  $L(t,V,\Delta)$ and $I_{\Delta}$ are given by Definition \ref{def-L} in section \ref{section-etimates-smooth-solu}.
\end{definition}

In addition we use the notations in Definition \ref{def-L0D0Q0} for $(u,v)$, and use the following for $(u^{\prime},v^{\prime})$,
\[
L_0^{\prime}(t,\Delta)=L(t,u^{\prime},\Delta)+L(t,v^{\prime},\Delta)
\]
and
\[
D_0^{\prime}(t,\Delta)=\int_{z_a(t)}^{b-t+t_0} |u^{\prime}(x,t)|^2|v^{\prime}(x,t)|^2 dx,
\]
\[ Q_0^{\prime}(t,\Delta)=
\int\int_{z_a(t)<x<y<b-t+t_0} |u^{\prime}(x,t)|^2|v^{\prime}(y,t)|^2 dxdy
\]
for $t\in I_{\Delta}$, and
\[r_0^{\prime}(x,t)=m(|u^{\prime}(x,t)|^2+|v^{\prime}(x,t)|^2)+8|\beta||u^{\prime}(x,t)|^2|v^{\prime}(x,t)|^2.\]
Moreover, (\ref{eq-dirac1}) and (\ref{eq-dirac2}) still hold for both $(u,v)$ and $(u^{\prime},v^{\prime})$, and  Lemmas in Section \ref{section-etimates-smooth-solu} also hold for  these two solution.

Now for any $T>0$, we can have the estimates on $F_1$ near the boundary $\Gamma_B$ as follows.
\begin{lemma}\label{lemma-bdry-BonyFunctional}
Suppose that $\Delta\subset R^1\times [0,T]$ and $\Delta\cap \Gamma_B\neq \emptyset$. Then, there exist constants  $\delta_0>0$, $K_1>0$ and $C_1>0$ such that if $L_0(t_0,\Delta)\le \delta_0$ and  $L_0^{\prime}(t_0,\Delta)\le \delta_0$ then there holds the following,
\begin{equation}\label{eq-ineqF1}
\frac{d}{dt}F_1(t,\Delta)\le -D_1(t,\Delta)+\big[O(1)+C_1\Lambda_1(t,\Delta)+C_1\Lambda_2(t)\big]F_1(t,\Delta)\end{equation}
for $t\in I_{\Delta}$ with  $z(t)\neq a+t-t_0$,
where
\[ \Lambda_1(t,\Delta)=4m\delta_0+8|\beta|(D_0(t,\Delta)+D^{\prime}_0(t,\Delta)),\]
and
\[\Lambda_2(t)=\left\{\begin{array}{ll} (1-z_t(t))|\lambda(t)|^2(|v(z(t),t)|^2 +|v^{\prime}(z(t),t)|^2), & \mbox{if $a+t-t_0\le z(t)$,} \\
0, & \mbox{if $a+t-t_0> z(t)$.}\end{array}
\right. \]
Here the constants $K_1>1$,  $\delta_0$ and $ C_1$ depend only on $T$.
\end{lemma}
{\it Proof.} It suffices to prove lemma for two cases.

Case 1: The boundary $\Gamma_B$ and the right edge of $\Delta$ intersec at some point $(z(\tau_1),\tau_1)$, see Figure \ref{fig-case1}. Then $I_{\Delta}= [t_0, \tau_1]$, and  $z_a(t)=z(t)$ for $t\in I_{\Delta}$.

For $t\in (t_0, \tau_1)$,  by Lemma \ref{lemma-nonlinearterm1} and Lemma \ref{lemma-nonlinearterm2},  we use (\ref{eq-dirac3}) and (\ref{eq-dirac2}) for both $(u,v)$ and $(u^{\prime},v^{\prime})$ to derive that
\begin{eqnarray*}
&\,& \frac{d}{dt} \int\int_{z(t)<x<y<b-t+t_0} |U(x,t)|^2(|v(y,t)|^2+|v^{\prime}(y,t)|^2) dxdy
\\ &\le & -\int\int_{z(t)<x<y<b-t+t_0} (|U(x,t)|^2)_x(|v(y,t)|^2+|v^{\prime}(y,t)|^2) dxdy \\ &\,& +\int\int_{z(t)<x<y<b-t+t_0} |U(x,t)|^2(|v(y,t)|^2+|v^{\prime}(y,t)|^2)_y dxdy \\ &\,& +\int\int_{z(t)<x<y<b-t+t_0}r_1(x,t)(|v(y,t)|^2+|v^{\prime}(y,t)|^2) dxdy \\ &\,&
+\int\int_{z(t)<x<y<b-t+t_0} |U(x,t)|^2(r_0(y,t)+r_0^{\prime}(y,t))dxdy \\ &\,&-\int_{z(t)}^{b-t+t_0}|U(x,t)|^2(|v(b-t+t_0,t)|^2+|v^{\prime}(b-t+t_0,t)|^2)dx \\&\,&-z_t(t)\int^{b-t+t_0}_{z(t)}|U(z(t),t)|^2(|v(y,t)|^2+|v^{\prime}(y,t)|^2)dy \\
 &\le&
-2\int^{b-t+t_0}_{z(t)}|U(y,t)|^2\big( |v(y,t)|^2+|v^{\prime}(y,t)|^2\big)dy \\ &\,& +(1-z_t(t))|\lambda(t)|^2|V(z(t),t)|^2 \int^{b-t+t_0}_{z(t)}(|v(y,t)|^2+|v^{\prime}(y,t)|^2)dy \\
&\,& +(mL_1(t)+c_*D_1(t))\int^{b-t+t_0}_{z(t)}(|v(y,t)|^2+|v^{\prime}(y,t)|^2)dy \\
&\,&+ (mL_0(t)+8|\beta|D_0(t)+mL_0^{\prime}(t)+8|\beta|D_0^{\prime}(t))\int_{z(t)}^{b-t+t_0}|U(x,t)|^2 dx,
\end{eqnarray*}
while  by Lemma \ref{lemma-nonlinearterm1} and Lemma \ref{lemma-nonlinearterm2},  we use (\ref{eq-dirac4}) and (\ref{eq-dirac1}) for both $(u,v)$ and $(u^{\prime},v^{\prime})$ to derive that
\begin{eqnarray*}
&\,& \frac{d}{dt} \int\int_{z(t)<x<y<b-t+t_0} (|u(x,t)|^2+|u^{\prime}(x,t)|^2)|V(y,t)|^2 dxdy \\
&\le& -\int\int_{z(t)<x<y<b-t+t_0}(|u(x,t)|^2+|u^{\prime}(x,t)|^2)_x|V(y,t)|^2 dxdy \\
&\,&
+\int\int_{z(t)<x<y<b-t+t_0} (|u(x,t)|^2+|u^{\prime}(x,t)|^2)(|V(y,t)|^2)_y dxdy \\
&\,& +\int\int_{z(t)<x<y<b-t+t_0} (r_0(x,t)+r_0^{\prime}(x,t))|V(y,t)|^2 dxdy \\
&\,& +\int\int_{z(t)<x<y<b-t+t_0} (|u(x,t)|^2+|u^{\prime}(x,t)|^2)r_1(y,t)dxdy \\
&\,& -z_t(t)\int^{b-t+t_0}_{z(t)}(|u(z(t),t)|^2+|u^{\prime}(z(t),t)|^2)|V(y,t)|^2 dy \\
&\,& -\int_{z(t)}^{b-t+t_0}(|u(x,t)|^2+|u^{\prime}(x,t)|^2)|V(b-t+t_0,t)|^2 dx \\
&\le& -2\int_{z(t)}^{b-t+t_0}(|u(x,t)|^2+|u^{\prime}(x,t)|^2|V(x,t)|^2)dx \\
&\,& +(1-z_t(t))(|u(z(t),t)|^2+|u^{\prime}(z(t),t)|^2)\int_{z(t)}^{b-t+t_0} |V(x,t)|^2dx \\
&\,&+(mL_1(t)+c_*D_1(t))\int_{z(t)}^{b-t+t_0}(|u(x,t)|^2+|u^{\prime}(x,t)|^2)dx \\
&\,&+ (mL_0(t)+8|\beta|D_0(t)+mL_0^{\prime}(t)+8|\beta|D_0^{\prime}(t))
\int_{z(t)}^{b-t+t_0}|V(x,t)|^2 dx.
\end{eqnarray*}

Collecting these two inequalities, we have the estimates on $Q_1$ as follows,
\begin{eqnarray*}
\frac{d}{dt}Q_1(t) &\le& \big[-2+c_*(L_0(t)+L_0^{\prime}(t)) \big]D_1(t)+[q_1(t)+q_2(t)]L_1(t) \\
&\,&+(1-z_t(t))|\lambda(t)|^2|V(z(t),t)|^2(L_0(t)+L^{\prime}_0(t)),
\end{eqnarray*}
where
\[ q_1(t)=2m(L_0(t)+L_0^{\prime}(t))+8|\beta|(D_0(t)+D_0^{\prime}(t)),\]
and
\[ q_2(t)=(1-z_t(t))|\lambda(t)|^2(|v(z(t),t)|^2+|v^{\prime}(z(t),t)|^2).\]

For the functional $L_1$, by (\ref{eq-dirac3}) and by Lemma \ref{lemma-nonlinearterm2}, we have
\begin{eqnarray*}
\frac{d}{dt}L(t,U) &\le & -|U(b-t+t_0,t)|^2-z_t(t)|U(z(t),t)|^2 \\ &\,& -\int^{b-t+t_0}_{z(t)}(|U(x,t)|^2)_x dx + \int^{b-t+t_0}_{z(t)} r_1(x,t)dx \\ &\le & (1-z_t(t))|\lambda(t)|^2|V(z(t),t)|^2 +\int^{b-t+t_0}_{z(t)} r_1(x,t)dx,
\end{eqnarray*}
while by (\ref{eq-dirac4}) and by Lemma \ref{lemma-nonlinearterm2}, we have
\begin{eqnarray*}
\frac{d}{dt}L(t,V) &\le& -|V(b-t+t_0,t)|^2-z_t(t)|V(z(t),t)|^2 \\ &\,& +\int_{z(t)}^{b-t+t_0}(|V(x,t)|^2)_xdx +\int_{z(t)}^{b-t+t_0} r_1(x,t)dx \\
&\le& -(1+z_t(t))|V(z(t),t)|^2+\int_{z(t)}^{b-t+t_0} r_1(x,t)dx.
\end{eqnarray*}

Then we have the following estimate on $L_1$,
\begin{eqnarray*}
\frac{d}{dt}\big( L(t,U)+K_1L(t,V)\big) &\le & \big[(1-z_t(t))|\lambda(t)|^2-K_1(1+z_t(t))\big]|V(z(t),t)|^2 \\ &\,&+(1+K_1)\int_{z(t)}^{b-t+t_0} r_1(x,t)dx \\
&\le&-2|V(z(t),t)|^2 +(1+K_1)\big( mL_1(t)+c_*D_1(t)\big).
\end{eqnarray*}
Here the constant $K_1>1$ is chosen so that
\[ \max_{t\in [0,T]} \big[(1-z_t(t))|\lambda(t)|^2-K_1(1+z_t(t))\big] <-2.\]

Now, with the above estimates on $Q_1$ and $L_1$, we use Lemma \ref{lemma-local-charge} to derive the following,
\begin{eqnarray*}
\frac{d}{dt}F_1(t) &\le & \big\{(1+K_1)c_*+ [-2+c_*(L_0(t)+L^{\prime}(t)) ]C_1\big\}D_1(t)  \\
&\,& +\big[ -2+C_1(1-z_t(t))|\lambda(t)|^2(L_0(t)+L^{\prime}_0(t))\big] |V(z(t),t)|^2 \\
&\,& +\big[(1+K_1)m +C_1(q_1(t)+q_2(t))\big]L_1(t)  \\
&\le& -D_1(t)+\big[(1+K_1)m +C_1(\Lambda_1(t)+\Lambda_2(t))\big]L_1(t)
\end{eqnarray*}
for $ L_0(t_0)\le \delta_0$ and $ L_0^{\prime}(t_0)\le \delta_0$,
where we choose $\delta_0>0$ and $C_1>0$ so that
\[ -2+2c_*\delta_0<-1, \quad (1+K_1)c_*-C_1<-1,\]
and
\[ -2+2C_1\delta_0\max_{0\le t\le T}(1-z_t(t))|\lambda(t)|^2<-1.\]
Therefore (\ref{eq-ineqF1}) is proved for Case 1.

Case 2: The boundary $\Gamma_B$ and the left edge of $\Delta$ intersect at $(z(\tau_2), \tau_2)$. Then, $I_{\Delta}=[t_0,t_0+\frac{b-a}{2}]$, and $z_a(t)=z(t)$ for $t_0\le t\le\tau_2$, $z_a(t)=a-t_0+t$ for $\tau_2\le t \le  t_0+\frac{b-a}{2}$.  The proof can be carried out in the same way as in Case 1 for $t\neq \tau_2$.
Thus the proof is complete.$\Box$

\begin{remark}
For the case that $\Delta\subset \Omega(T)$, we have similar estimates on $F_1$ without boundary terms, see \cite{Zhang-zhao} for the proof, where only $D_1(t,\Delta)$ makes contribution to the control on $F_1$. For the case that $\Delta\cap \Gamma_B\neq\emptyset$,  both $Q_1(t,\Delta)$ and $L(t,V,\Delta)$ are needed to give the control on $F_1$.
\end{remark}

As conclusion of the above argument, we get the stability result for smooth solutions for any $T>0$.
\begin{proposition}\label{lemma-stability} Suppose that $\Delta\subset R^1\times [0,T]$ with $b>z(t_0)$, and suppose that $L_0(t_0,\Delta)\le \delta_0$,  $L_0^{\prime}(t_0,\Delta)\le \delta_0$. Then
for $t\in I_{\Delta}$, there holds the following
\begin{eqnarray*} &\,& \int_{z_a(t)}^{b+t_0-t}(|u(x,t)-u^{\prime}(x,t)|^2+|v(x,t)-v^{\prime}(x,t)|^2)dx \\&\le& C_4 \int_{\max\{z(t_0),a\}}^b(|u(x,t_0)-u^{\prime}(x,t_0)|^2+|v(x,t_0)-v^{\prime}(x,t_0)|^2)dx,\end{eqnarray*}
and
\begin{eqnarray*}
&\,& \int\int_{\Delta\cap\Omega}(|uv-u^{\prime}v^{\prime}|^2)dxdt \\&\le& C_4 \int_{\max\{z(t_0),a\}}^b(|u(x,t_0)-u^{\prime}(x,t_0)|^2+|v(x,t_0)-v^{\prime}(x,t_0)|^2)dx,\end{eqnarray*}
\begin{eqnarray*}
&\,& \int\int_{\Delta\cap\Omega}(|u\overline{v}-u^{\prime}\overline{v^{\prime}}|^2)dxdt \\&\le& C_4 \int_{\max\{z(t_0),a\}}^b(|u(x,t_0)-u^{\prime}(x,t_0)|^2+|v(x,t_0)-v^{\prime}(x,t_0)|^2)dx.\end{eqnarray*}
Here the constant $C_4$ depends only on $T$ and $E_0$.
\end{proposition}
{\it Proof.} It suffices to prove lemma for two cases.

Case 1: $a<z(t_0)<b$, that is, $\Delta \cap \Gamma_B\neq\emptyset$. Then taking the integral of (\ref{eq-ineqF0}) in Lemma \ref{lemma-bdry-Q} over $I_{\Delta}$, we have
\begin{eqnarray*} & \int_{I_{\Delta}} \big(D_0(t,\Delta)+D^{\prime}_0(t,\Delta) +|v(z_a(t),t)|^2+|v^{\prime}(z_a(t),t)|^2\big)dt
 \\ &\le O(1) (L_0(t_0,\Delta)+L_0^{\prime}(t_0,\Delta))+O(1)\delta_0 T,
\end{eqnarray*}
 which leads to
 \[ \int_{I_{\Delta}} \big(\Lambda_1(t,\Delta)+\Lambda_2(t,\Delta) \big)dt \le C(T) \]
 for some constant $C(T)>0$ depending on $T$.

 Therefore, we use Lemma \ref{lemma-bdry-BonyFunctional} to deduce that
 \begin{eqnarray*} F_1(t)&\le& \exp(\int_{I_{\Delta}}[O(1)+C_1\Lambda_1(s,\Delta)+C_1\Lambda_2(s,\Delta)]ds)F_1(t_0)
 \\ &\le& C^{\prime}(T)
 \int_{z(t_0)}^b(|u(x,t_0)-u^{\prime}(x,t_0)|^2+|v(x,t_0)-v^{\prime}(x,t_0)|^2)dx \end{eqnarray*}
 for $t\in I_{\Delta}$,
 and
 \begin{eqnarray*}
 \int_{I_{\Delta}} D_1(t,\Delta)dt &\le& F_1(t_0) \\ &+&(\int_{I_{\Delta}}[O(1)+C_1\Lambda_1(s,\Delta)+C_1\Lambda_2(s,\Delta)]ds)T\max_{t_0\le t\le t_0+\frac{b-a}{2}}F_1(t)
 \\&\le& F_1(t_0)+C^{\prime\prime}(T) \max_{t_0\le t\le t_0+\frac{b-a}{2}}F_1(t),
 \end{eqnarray*}
 which lead to the result for Case 1. Here the constants $C^{\prime}(T)$ and $C^{\prime\prime}(T)$ depend only on $T$.

 Case 2: $z(t_0)<a$, that is, $\Delta \cap \Gamma_B=\emptyset$ and $\Delta\subset\Omega(T)$. Then $z_a(t)=a+t-t_0$. The result for this case has been proved in \cite{Zhang-zhao}, and its proof can be carried out in the same way as above. Therefore the proof is complete.
 $\Box$

\section{Convergence of global classical solutions}\label{section-convergence}

Choose a sequence of smooth functions
\[ (u^{(k)}_0,v^{(k)}_0)\in C^{\infty}_c(0,\infty), \quad k=1,2,\cdot,\]
such that \[(u^{(k)}_0,v^{(k)}_0)\to (u_0,v_0) \quad in \, L^2_{loc}(0,\infty)\] as $m\to \infty$. Theorem \ref{thm-1} implies that there is a sequence of classical solutions, $(u^{(k)},v^{(k)})\in C^1(\overline{\Omega})$,  $k=1,2,\cdots,$ to (\ref{eq-dirac}), which satisfy boundary condition (\ref{eq-dirac-bc}) and take $(u^{(k)}_0,v^{(k)}_0)$ as their initial data respectively. And $supp(u^{(k)}(\cdot, t),v^{(k)}(\cdot, t))$ has bounded support in $R^1$ for any $t\ge 0$ and $k\ge 0$.

We consider the convergence of $\{(u^{(k)},v^{(k)})\}_{k=0}^{\infty}$  on $\Delta(-A,A,0)\cap\Omega$ for any $A>0$. To this end, we first give the estimate on $L^2$ norm of solution over small interval $[a,b]\cap [z(t), A-t]$ for any $a$ and $b$.
\begin{lemma}\label{lemma-uniformconti} There is a constant $r>0$ such that if $0<b-a\le 4r$ and $b\le A$ then
\[ \sup_{k\ge 1} \int_{\max\{z(t),a\}}^{\min\{b,A-t\}} (|u^{(k)}(x,t)|^2+|v^{(k)}(x,t)|^2) dx\le \delta_0\]
for $t\in [0,A]$ with $z(t)\le b$.
\end{lemma}
{\it Proof.}It is obvious that
\[\lim_{k\to\infty} \int_{0}^A(|u^{(k)}_0-u_0|^2+|v^{(k)}_0-v_0|^2)dx=0.\]
As in \cite{Zhang-zhao}, we choose $r>0$ such that
\[ \exp(mA+8|\beta|j_0)\big( \int_{\max\{z(0),a\}}^b (|u^{(k)}_0(x)|^2+|v^{(k)}_0(x)|^2) dx +mj_0(b-a)\big) \le \frac{\delta_0}{8} \]
and
\[\exp(2mA+16|\beta|j_0)\big(\int_{\max\{z(0),a\}}^b |v_0^{(k)}(2p(x)+x)|^2dx+mj_0(b-a)\big) \le \frac{\delta_0}{8} \]
for  $|b-a|\le 4r$ and for $k=0,1,2\cdots$. Here for  simplification $(u^{(0)}_0,v^{(0)}_0)=(u_0,v_0)$.

Then with the pointwise estimates along the characteristics in Lemma \ref{lemma-pointwise-v} and Lemma \ref{lemma-pointwise-u}, we can deduce the desired result. The proof is complete. $\Box$

Now application of Proposition \ref{lemma-stability} and Lemma \ref{lemma-uniformconti} to any pair of smooth solutions $(u^{(k)},v^{(k)})$ and $(u^{(n)},v^{(n)})$ gives the following.
\begin{lemma}\label{lemma-interiorestimate}
Suppose that $\Delta (a,b,\tau)\subset\Delta(-A,A,0)$ with $0<b-a\le 4r$ and $\Delta(a,b,\tau)\cap \Omega\neq \emptyset$. Then there exists a constant $C(A)>0$ such that
\begin{eqnarray*} &\,& \int\int_{\Delta(a,b,\tau)\cap\Omega}(|u^{(k)}-u^{(n)}|^2+|v^{(k)}-v^{(n)}|^2)dxdt \\
&+&\int\int_{\Delta(a,b,\tau)\cap\Omega}(|u^{(k)}v^{(k)}-u^{(n)}v^{(n)}|^2
+|u^{(k)}\overline{v^{(k)}}-u^{(n)}\overline{v^{(n)}}|^2)dxdt \\
&\le& C(A) \int_{\max(z(0),a)}^b(|u^{(k)}(x,\tau)-u^{(n)}(x,\tau)|^2+|v^{(k)}(x,\tau)-v^{(n)}(x,\tau)|^2)dx
\end{eqnarray*}
for any $k\ge 1$ and $n\ge 1$.
Here the constant $C(A)$ depends only on $A$ and $E_0$; the constant $r>0$ is given by Lemma \ref{lemma-uniformconti}.
\end{lemma}

In the next, we prove the convergence of $\{(u^{(k)},v^{(k)})\}^{\infty}_{k=0}$ on $\Delta(-A,A,0)\cap\Omega$  by the induction step as follows.

Denote
\[ \Omega(A,\tau)=\Delta(-A,A,0)\cap\Omega\cap\{(x,t)\,|\, 0\le t\le \tau\}, \quad 0\le \tau\le A.\]
\begin{figure}[ht]
\begin{center}
\unitlength=10mm
\begin{picture}(10,3.5)
\thicklines
\put(0,0){\line(1,0){10}}
\put(1,0){\line(5,4){4}}
\put(9,0){\line(-5,4){4}}
\put(0,-0.5){$(-A,0)$}
\put(8.5,-0.5){$(A,0)$}\put(10.3,0){$t=0$}\put(0,1){\line(1,0){10}}\put(0.5,1.2){$(\tau-A,\tau)$}\put(7.8,1.2){$(A-\tau,\tau)$}
\qbezier(5,-0.5)(3.5,3)(3.5,4)
\put(4.5,1.3){$t=\tau$}\put(5,0.5){$\Omega(A,\tau)$}
\end{picture}
\caption{Domain $\Omega(A,\tau)$}\label{fig-induction step}
\end{center}
\end{figure}
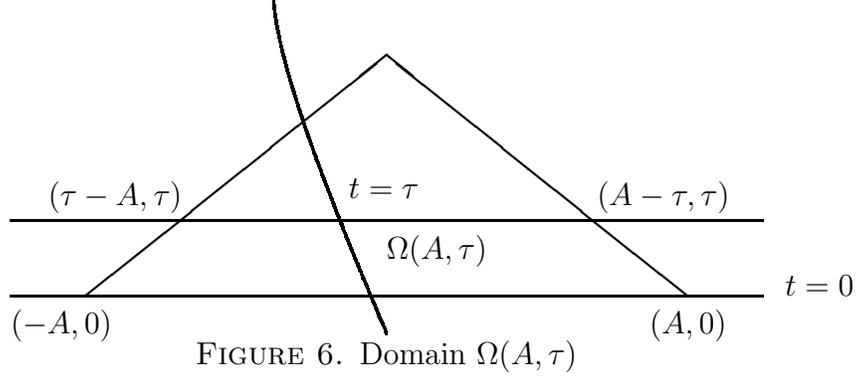

\begin{lemma}\label{lemma-induction}
Suppose that
\[\lim_{m,n\to\infty}\int\int_{\Omega(A,\tau)}(|u^{(m)}-u^{(n)}|^2+|v^{(m)}-v^{(n)}|^2)dxdt=0,\]
\[ \lim_{m,n\to\infty}\int\int_{\Omega(A,\tau)}(|u^{(m)}v^{(m)}-u^{(n)}v^{(n)}|^2)dxdt=0\]
and
\[ \lim_{m,n\to\infty}\int\int_{\Omega(A,\tau)}(|u^{(m)}\overline{v^{(m)}}-u^{(n)}\overline{v^{(n)}}|^2)dxdt=0\]
for $\tau\in [0,A-r]$.
Then
\[\lim_{m,n\to\infty}\int\int_{\Omega(A,\tau+r)}(|u^{(m)}-u^{(n)}|^2+|v^{(m)}-v^{(n)}|^2)dxdt=0,\]
\[ \lim_{m,n\to\infty}\int\int_{\Omega(A,\tau+r)}(|u^{(m)}v^{(m)}-u^{(n)}v^{(n)}|^2)dxdt=0\]
and
\[ \lim_{m,n\to\infty}\int\int_{\Omega(A,\tau+r)}(|u^{(m)}\overline{v^{(m)}}-u^{(n)}\overline{v^{(n)}}|^2)dxdt=0.\]
Here $r$ is given in Lemma \ref{lemma-uniformconti}.
\end{lemma}
{\it Proof.}
We choose a finite number of subintervals, $[a_j,b_j]$, $j=1,2,\cdots,J$, with $[a_j,b_j]\cap [\max\{\tau-A,z(\tau)\},A-\tau]\neq\emptyset$ and $b_j-a_j=4r$, such that
\[ \overline{\Omega(A,\tau+r)\backslash\Omega(A,\tau)}\subset \cup_{j=1}^{J}\Delta(a_j,b_j,\tau^{\prime}),\]
where $\tau^{\prime}=\tau-\frac{r}{8}$, and $J\le \frac{4A}{r}+1$.

For $1\le j\le J$, by Proposition \ref{lemma-stability} and Lemma \ref{lemma-uniformconti}, we have
\[\lim_{m,n\to\infty}\int\int_{\Delta(a_j,b_j,\tau^{\prime})\cap\Omega}(|u^{(m)}-u^{(n)}|^2+|v^{(m)}-v^{(n)}|^2)dxdt=0,\]
\[ \lim_{m,n\to\infty}\int\int_{\Delta(a_j,b_j,\tau^{\prime})\cap\Omega}(|u^{(m)}v^{(m)}-u^{(n)}v^{(n)}|^2)dxdt=0\]
and
\[ \lim_{m,n\to\infty}\int\int_{\Delta(a_j,b_j,\tau^{\prime})\cap\Omega}(|u^{(m)}\overline{v^{(m)}}-u^{(n)}\overline{v^{(n)}}|^2)dxdt=0.\]
Therefore we have the convergence of the sequences $\{ (u^{(m)},v^{(m)})\}_{m=1}^{\infty}$, $\{ u^{(m)}v^{(m)}\}_{m=1}^{\infty}$ and $\{ u^{(m)}\overline{v^{(m)}}\}_{m=1}^{\infty}$ in $L^2(\Omega(A,\tau+r))$ respectively. The proof is complete. $\Box$

Now we have the following convergence result.
\begin{proposition}\label{prop-convergence}
There exists a $(u,v)\in L^2_{loc}(\Omega)$ such that
\[\lim_{m\to\infty}|| (u^{(m)},v^{(m)})-(u,v)||_{L^2(\Delta(-A, A,0)\cap\Omega)}=0,\] and
\[\lim_{m\to\infty}\big(|| u^{(m)}v^{(m)}-uv||_{L^2(\Delta(-A, A,0)\cap\Omega)}+|| u^{(m)}\overline{v^{(m)}}-u\overline{v}||_{L^2(\Delta(-A, A,0)\cap\Omega)}\big)=0\] for any $A>0$.
\end{proposition}
{\it Proof.}
With the induction steps given by Lemma \ref{lemma-induction}, we have
\[\lim_{m,n\to\infty}\int\int_{(\Delta(-A, A,0)\cap\Omega)}(|u^{(m)}-u^{(n)}|^2+|v^{(m)}-v^{(n)}|^2)dxdt=0,\]
\[ \lim_{m,n\to\infty}\int\int_{(\Delta(-A, A,0)\cap\Omega)}(|u^{(m)}v^{(m)}-u^{(n)}v^{(n)}|^2)dxdt=0\]
and
\[ \lim_{m,n\to\infty}\int\int_{(\Delta(-A, A,0)\cap\Omega)}(|u^{(m)}\overline{v^{(m)}}-u^{(n)}\overline{v^{(n)}}|^2)dxdt=0,\]
for any $A>0$. These lead to the desired result. The proof is complete.$\Box$

\section{Proof of main results on strong solutions}\label{section-proof}

In the same way as in the proof of Lemma \ref{lemma-induction} and Proposition \ref{prop-convergence}, we can prove the following.

\begin{proposition}\label{prop-unique} Suppose that $\{u^{(m)}_j, v^{(m)}_j\}_{m=1}^{\infty}$, $j=1,2$, are two sequences of classical solution to (\ref{eq-dirac}) satisfy boundary condition (\ref{eq-dirac-bc}) with the following,
\[ \lim_{m\to\infty} \int_0^M(|u^{(m)}_1(x,0)-u^{(m)}_2(x,0)|^2+|v^{(m)}_1(x,0)-v^{(m)}_2(x,0)|^2)dx=0\] for some $M>0$.
Then,
\[\lim_{m\to\infty} \int\int_{\Delta(-M,M,0)\cap\Omega} (|u^{(m)}_1-u^{(m)}_2|^2+|v^{(m)}_1-v^{(m)}_2|^2)dxdt=0.\]
\end{proposition}

{\bf Proof of Theorem\ref{thm-3}.} The existence of solution $(u,v)$ is proved by Proposition \ref{prop-convergence}. Moreover, $(u,v)$ satisfies (\ref{eq-weaksolu-1}) and (\ref{lemma-pointwise-v}).

To prove the uniqueness, let $(u_j,v_j)$, $j=1,2$, be two strong solutions to (\ref{eq-dirac}-\ref{eq-dirac-bc}), and let $(u^{(m)}_j,v^{(m)}_j)$, $j=1,2$ be two sequences of classical solutions to (\ref{eq-dirac}) with boundary condition (\ref{eq-dirac-bc}), which are convergent to $(u_j,v_j)$, $j=1,2$, respectively in $L^2_{loc}(\Omega)$. Moreover,  the initial data $(u^{(m)}_j(x,0), v^{(m)}_j(x,0))$ are assumed to be convergent to $(u_0,v_0)$ for $j=1,2$.

Then by Proposition \ref{prop-unique}, we have
\[ \lim_{m\to\infty} \int\int_{\Delta(-A,A,0)\cap\Omega}(|u^{(m)}_1(x,0)-u^{(m)}_2(x,0)|^2+|v^{(m)}_1(x,0)-v^{(m)}_2(x,0)|^2)dx=0,\]
which yields that
\[ (u_1,v_1)(x,t)=(u_2,v_2)(x,t), \quad a.e. \, \, (x,t)\in \Delta(-A,A,0)\cap\Omega.\]
This leads to the uniqueness of the strong solution. The proof is complete.$\Box$

{\bf Proof of Theorem \ref{thm-4}.} Indeed the results hold for the classical solutions. Then by taking the limit, we can prove the result still hold for the strong solution. The proof is complete. $\Box$

\section*{ Acknowledgement}

This work was partially  supported by NSFC Project 11421061 and by the 111 Project
B08018.

\end{document}